\documentclass{article}
\usepackage{amssymb,amsmath,amsthm,graphicx,color,tikz}

\newtheorem{theorem}{Theorem}
\newtheorem{lemma}[theorem]{Lemma}
\newtheorem{proposition}[theorem]{Proposition}

\theoremstyle{definition}
\newtheorem{example}{Example}
\newtheorem{definition}{Definition}
\newtheorem{remark}{Remark}

\definecolor{myred}{RGB}{0, 0, 0}
\definecolor{mygreen}{RGB}{0, 0, 0}

\newcommand{\bx}[2]{{#1}_{#2}}

\newcommand{\btri}[3]{%
  \bx{#1}{\bssia} 
  \bx{#2}{\bssib} 
  \bx{#3}{\bssic}
}

\newcommand{\btrI}[3]{%
  \bx{#1}{\bssIc} 
  \bx{#2}{\bssIb} 
  \bx{#3}{\bssIa}%
}

\newcommand{\bssia}{{x,y}}
\newcommand{\bssib}{{x \utr y, z \otr y}}
\newcommand{\bssic}{{y,z}}

\newcommand{\bssIc}{{y \otr x, z \otr x}}
\newcommand{\bssIb}{{x,z}}
\newcommand{\bssIa}{{x \utr z, y \utr z}}
\newcommand{\pheq}{\phantom{=}\,}

\newcommand{\tri}[3]{%
  \pbx[\acol{a}]{#1}{\ssia} \cdot
  \pbx[\bcol{b}]{#2}{\ssib} \cdot
  \pbx[\ccol{c}]{#3}{\ssic}
}

\newcommand{\trI}[3]{%
  \pbx[\ccol{c}]{#1}{\ssIc} \cdot
  \pbx[\bcol{b}]{#2}{\ssIb} \cdot
  \pbx[\acol{a}]{#3}{\ssIa}%
}

\newcommand{\pbA}[2][0]{{A}_{{#1}, #2}}
\newcommand{\pbB}[2][0]{{B}_{{#1}, #2}}
\newcommand{\pbx}[3][0]{{#2}_{{#1}, #3}}
\newcommand{\pbxNoSubs}[2][0]{{#2}_{{#1}}}
\newcommand{\Ae}{\pbxNoSubs[0]{A}}
\newcommand{\Be}{\pbxNoSubs[0]{B}}
\newcommand{\Ao}{\pbxNoSubs[1]{A}}
\newcommand{\Bo}{\pbxNoSubs[1]{B}}

\newcommand{\ssia}{{x,y}}
\newcommand{\ssib}{{x \acol{\utr^a} y, z \ccol{\otr^c} y}}
\newcommand{\ssic}{{y,z}}

\newcommand{\ssIc}{{y \acol{\otr^a} x, z \bcol{\otr^b} x}}
\newcommand{\ssIb}{{x,z}}
\newcommand{\ssIa}{{x \bcol{\utr^b} z, y \ccol{\utr^c} z}}

\usetikzlibrary{
  knots,
  hobby,
  decorations.markings,
  decorations.pathreplacing,
  shapes.geometric,
  calc
}
\tikzset{
  show curve controls/.style={
    postaction=decorate,
    decoration={show path construction,
      curveto code={
        \draw [blue, dashed]
        (\tikzinputsegmentfirst) -- (\tikzinputsegmentsupporta)
        node [at end, draw, solid, red, inner sep=2pt]{};
        \draw [blue, dashed]
        (\tikzinputsegmentsupportb) -- (\tikzinputsegmentlast)
        node [at start, draw, solid, red, inner sep=2pt]{}
        node [at end, fill, blue, ellipse, inner sep=2pt]{}
        ;
      }
    }
  },
 show curve endpoints/.style={
    postaction=decorate,
    decoration={show path construction,
      curveto code={
        \node [fill, blue, ellipse, inner sep=2pt] at (\tikzinputsegmentlast) {};
      }
    }
  },
  add arrow/.style={postaction={decorate}, decoration={
      markings,
      mark=at position 0.2 with {\arrow{<}},
      mark=at position 0.6 with {\arrow{<}}}}
}

\newcommand{\acol}[1]{{\color{myred} #1}}
\newcommand{\bcol}[1]{{\color{mygreen} #1}}
\newcommand{\ccol}[1]{{\color{myred} #1}}

\DeclareMathOperator{\otr}{\overline{\triangleright}}
\DeclareMathOperator{\utr}{\underline{\triangleright}}
\def\overop{\otr}
\def\underop{\utr}

\begin{document}

\title{\Large \textbf{Kaestner Brackets}}

\author{Forest Kobayashi\footnote{Email: fkobayashi@hmc.edu. Supported by the Giovanni Borrelli Mathematics Fellowship.} \and
Sam Nelson \footnote{Email: Sam.Nelson@cmc.edu. Partially supported by
Simons Foundation collaboration grant $\#316709$.}
}
\date{}

  \maketitle

  \begin{abstract}
    We introduce \textit{Kaestner brackets}, a generalization of
    biquandle brackets to the case of parity biquandles. This infinite
    set of quantum enhancements of the biquandle counting invariant
    for oriented virtual knots and links includes the classical
    quantum invariants, the quandle and biquandle $2$-cocycle
    invariants and the classical biquandle brackets as special cases,
    coinciding with them for oriented classical knots and links but
    defining generally stronger invariants for oriented virtual knots
    and links. We provide examples to illustrate the computation of
    the new invariant and to show that it is stronger than the
    classical biquandle bracket invariant for virtual knots.
  \end{abstract}

\bigskip

  \parbox{5.5in}{

    \smallskip

  }

\section{Introduction}

While introducing virtual knots in \cite{K}, Kauffman noted that crossings in
Gauss code representations of knot diagrams have a well-defined
\textit{parity} which does not change under Reidemeister moves.
In \cite{KK}, this property of crossing parity was exploited to define
\textit{parity biquandles}, algebraic structures which may be understood as
biquandles with different biquandle operations at the even and odd parity 
crossings.

In \cite{NOR}, the second listed author and collaborators introduced
\textit{biquandle brackets}, skein invariants for biquandle-colored oriented
knots and links. The resulting infinite family of invariants includes the
classical quantum invariants as well as the quandle and biquandle cocycle
invariants as special cases, as well as new invariants. In \cite{NO}, the
second listed author and collaborator introduced a graphical calculus for
computing biquandle bracket invariants using \textit{trace diagrams} to allow
for recursive computation as opposed to state-sum style computation of
biquandle bracket invariants.

Biquandle bracket invariants are well-defined for oriented virtual
knots and links by the usual convention of ignoring virtual crossings
when computing biquandle colorings and states in the state-sum
expansion of the bracket (or, equivalently, working over a supporting
surface of minimal genus).  In this paper we
apply the biquandle bracket idea to the case of parity
biquandles, obtaining an algebraic structure 
 with associated virtual link invariants 
  that we call \textit{Kaestner brackets}. This infinite family of
oriented virtual link invariants includes classical
biquandle bracket invariants as a special case --
indeed, coinciding with them for classical knots -- but includes
new, generally stronger invariants for virtual knots and links.

The paper is organized as follows. In Section \ref{BG} we review the
basics of biquandles and parity biquandles, and we recall a few facts
about signed Gauss codes. In Section \ref{BB} we review the basics of
biquandle brackets. In Section \ref{KB} we generalize biquandle
brackets to the case of parity biquandles, defining Kaestner brackets.
We provide examples of the computation of a \textit{Kaestner bracket
polynomial} and collect some results from 
\texttt{python} computations. In particular, we show that Kaestner bracket
invariants can distinguish virtual knots which are not distinguished
by the corresponding classical biquandle bracket. We conclude in
Section \ref{Q} with some questions for future research.

\section{Biquandles and Parity Biquandles}\label{BG}

We begin with a definition; see \cite{EN} and the references therein for more.

\begin{definition}
Let $X$ be a set. A \textit{biquandle structure} on $X$ consists of binary
operations $\utr$, $\otr : X \times X \to X$ such that
\begin{itemize}
\item[(i)] For all $x \in X$, $x \utr x = x \otr x$,
\item[(ii)] For all $x,y \in X$, the maps $\alpha_y, \beta_y : X \to X$
and $S : X \times X \to X \times X$ defined by
\[\alpha_y(x)=x\otr y, \beta_y(x)=x\utr y\quad \mathrm{and}\quad
S(x,y)=(y\otr x,x\utr y) \]
are all invertible, and
\item[(iii)] For all $x,y,z \in X$, we have the following
\textit{exchange laws:}
\begin{eqnarray*}
(x \utr y) \utr (z \utr y) &= & (x \utr z) \utr (y \otr z) \\
(x \utr y) \otr (z \utr y) &= & (x \otr z) \utr (y \otr z) \\
(x \otr y) \otr (z \otr y) &= & (x \otr z) \otr (y \utr z).
\end{eqnarray*}
\end{itemize}
\end{definition}

\begin{definition}
Let $L$ be an oriented virtual link diagram and let $X$ be a biquandle. An
\textit{$X$-coloring} or \textit{biquandle coloring of $L$ by $X$}
is a labeling of all of the semiarcs of $L$ with elements of $X$
as shown:
\[\includegraphics{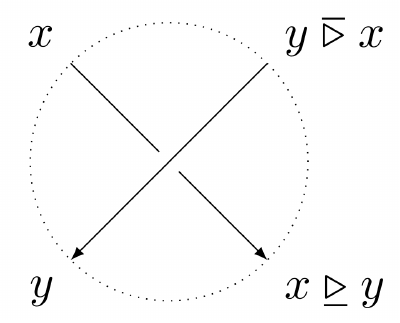}\quad \includegraphics{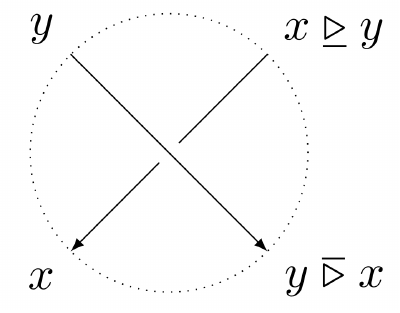}\]
Note that at virtual crossings the biquandle colors do not change, i.e., 
we have
\[\includegraphics{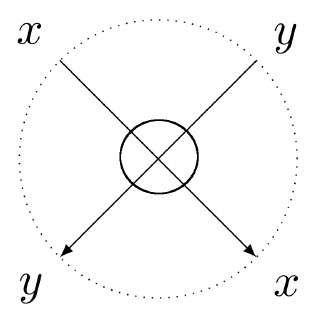}.\]
\end{definition}

The biquandle axioms are chosen so that given a biquandle coloring
of an oriented virtual link diagram before a Reidemeister move, there is a
\textit{unique} biquandle coloring of the diagram after the move
which agrees with the original coloring outside the neighborhood
of the move. In particular, we have the following result:

\begin{theorem}
Let $X$ be a finite biquandle, and let $L$ be an oriented virtual link diagram.
Then the number of $X$-colorings of an oriented virtual link is an
integer-valued invariant.
\end{theorem}

\begin{example}
Let $R$ be a commutative ring with identity and two distinguished units
$t,s$. Then $R$ is a biquandle (called an \textit{Alexander biquandle})
under the operations
\[ x \utr y = tx + (s - t)y \quad\quad x \otr y = s x. \]
\end{example}

\begin{example} In fact, we note that we can relax the ``commutative ring
with identity'' condition somewhat; any module over a ring $R$ with
commuting units $t,s$ forms an Alexander biquandle. For example, the ring
$M_n(\mathbb{F})$
of $n\times n$ matrices over a field $\mathbb{F}$ becomes an Alexander
biquandle with a choice of two commuting invertible matrices $t,s$.
When $\mathbb{F}$ is finite, so is the resulting biquandle.
\end{example}

\begin{example}
Let $X=\{1,2,\dots,n\}$. We can specify a biquandle structure on $X$ by
giving the operation tables of $\utr$ and $\otr$ explicitly, e.g.
\[\begin{array}{c|ccc}
\utr & 1 & 2 & 3 \\ \hline
1 & 1 & 3 & 2\\
2 & 3 & 2 & 1 \\
3 & 2 & 1 & 3 \\
\end{array}
\quad
\begin{array}{c|ccc}
\otr & 1 & 2 & 3 \\ \hline
1 & 1 & 1 & 1\\
2 & 2 & 2 & 2 \\
3 & 3 & 3 & 3
\end{array}.
\]
\end{example}

Next, we recall \emph{parity biquandles} (see  \cite{KK}):
  \begin{definition}
    Let $X$ be a set together with four binary operations $\underop^0,
    \overop^0, \underop^1, \overop^1$. Then we call $(X, \underop^0, \overop^0,
    \underop^1, \overop^1)$ a \emph{parity biquandle} iff
    \begin{itemize}
      \item $(X, \underop^0, \overop^0)$ is a biquandle (note, $(X, \underop^1,
        \overop^1)$ need not be a biquandle),
      \item For all $x,y \in X$, the maps $\alpha_y^1, \beta_y^1 : X \to X$
and $S^1 : X \times X \to X \times X$ defined by
\[\alpha_y^1(x)=x\otr^1 y, \beta_y^1(x)=x\utr^1 y\quad \mathrm{and}\quad
S^1(x,y)=(y\otr^1 x,\ x\utr^1 y) \]
are all invertible, and
      \item For all $(\acol{a},\bcol{b},\ccol{c}) \in
        \{(\acol{1},\bcol{1},\ccol{0}), (\acol{1},\bcol{0},\ccol{1}),
          (\acol{0},\bcol{1},\ccol{1})\}$ and for all $x,y,z \in X$, we have the
        \emph{mixed exchange laws}:
        \begin{align*} \textstyle
          (z \acol{\overop^a} y) \bcol{\overop^b} (x \ccol{\overop^c} y)
          &=  \textstyle
            (z \bcol{\overop^b} x) \acol{\overop^a} (y \ccol{\underop^c} x) \\
          \textstyle
          (x \acol{\overop^a} y) \bcol{\underop^b} (z \ccol{\overop^c} y)
          &=  \textstyle
            (x \bcol{\underop^b} z) \acol{\overop^a} (y \ccol{\underop^c} z) \\
          \textstyle
          (y \acol{\underop^a} x) \bcol{\underop^b} (z \ccol{\overop^c} x)
          &=  \textstyle
            (y \bcol{\underop^b} z) \acol{\underop^a} (x \ccol{\underop^c} z)
        \end{align*}
    \end{itemize}
  \end{definition}

\begin{example}
Biquandles are parity biquandles where $\underop^0 =  \underop^1$
and $\overop^0 = \overop^1$.
\end{example}

\begin{example}\label{ex:pb1}
The operation tables
\[
\begin{array}{r|rrr}
\utr^0 & 1 & 2 & 3 \\ \hline
1 & 3 & 1 & 3\\
2 & 2 & 2 & 2\\
3 & 1 & 3 & 1
\end{array}\quad
\begin{array}{r|rrr}
\otr^0 & 1 & 2 & 3 \\ \hline
1 & 3 & 1 & 3\\
2 & 2 & 2 & 2\\
3 & 1 & 3 & 1
\end{array}\quad
\begin{array}{r|rrr}
\utr^1 & 1 & 2 & 3 \\ \hline
1 & 3 & 1 & 3\\
2 & 2 & 2 & 2\\
3 & 1 & 3 & 1
\end{array}\quad
\begin{array}{r|rrr}
\otr^1 & 1 & 2 & 3 \\ \hline
1 & 1 & 3 & 1\\
2 & 2 & 2 & 2\\
3 & 3 & 1 & 3
\end{array}
\]
define a parity biquandle structure on the set $\{1,2,3\}$.
\end{example}

We note the following results; see \cite{KK}.

  \begin{lemma}
    Let $L$ be an oriented virtual knot or link diagram and let $c_i$ be a
    crossing in $L$.
    Then virtual Reidemeister moves preserve the parity of $c_i$.
  \end{lemma}

  \begin{lemma}
    Let $L$ be as above. Then for any Reidemeister III move, either all three
    crossings are even, or two are odd and one is even.
  \end{lemma}

Next, we briefly recall \textit{signed Gauss codes}; see \cite{K,GPV}.

\begin{definition}
Let $L$ be an oriented link represented as a planar diagram with $n$ crossings. Then we encode $L$ in a string of symbols by the following scheme:
    \begin{itemize}
      \item Pick some starting point $p_0$ on $K$ and begin transversing $K$ in the positive direction as specified by the orientation. Label new crossings as
        with $1, \ldots, n$ in the order that they're encountered. Each crossing
        should be visited exactly twice; we only label a crossing the first
        time.
      \item Whenever we encounter a crossing, we record three pieces of
        information:
        \begin{itemize}
          \item The crossing label,
          \item whether we're on the under/overstrand, and
          \item the sign of the crossing.
        \end{itemize}
        We can write this compactly by $k_{x}^\epsilon$, where $k \in \{1,
          \ldots, n\}$ is the label we've assigned our crossing, $x \in \{u,
          o\}$ denotes whether we're on the \textbf{u}nderstrand or
        \textbf{o}verstrand, and $\epsilon \in \{+, -\}$ denotes the sign of
        $k$.
    \end{itemize}
Repeat for each component of $L$. The resulting set of ordered strings is
the \emph{signed Gauss code} of $L$.
  \end{definition}
  \begin{example}
    \label{ex:GaussCode}
    The signed Gauss code
\[1^-_u, 2^-_o, 3^-_u, 4^-_o, 5^-_u, 6^-_o, 7^-_u, 1^-_o, 6^-_u, 5^-_o, 4^-_u, 3^-_o, 2^-_u, 7^-_o\]
corresponds to a diagram for the knot $7_2$:
\[ \scalebox{1.2}{\includegraphics{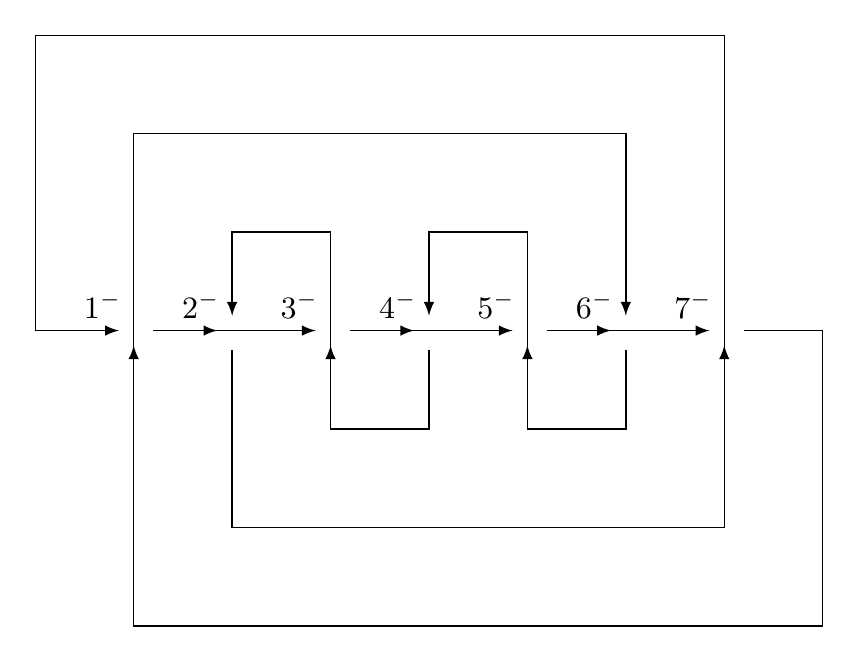}}\]
  \end{example}

\begin{definition}
Let $G$ be a signed Gauss code. A crossing in $G$ has \textit{even parity} if 
the number of crossing labels between its over and under instances is even and
\textit{odd parity}  if the number of crossing labels between the its over and 
under instances is odd.
\end{definition}

  Recall (see \cite{GPV,K}) that the Reidemeister moves on signed Gauss codes
  are as follows:
    \begin{itemize}
      \item Reidemeister I moves correspond to insertion/deletion of an adjacent
        pair: \label{ax:reidGaussCodeI}
        \[
          \scalebox{0.7}{\includegraphics{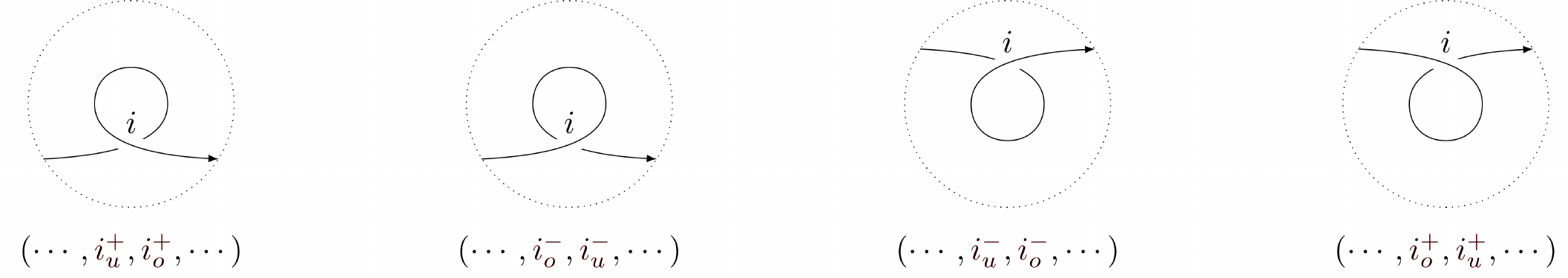}}
        \]
      \item Reidemeister II moves correspond to insertion/deletion of a pair of
        pairs: \label{ax:reidGaussCodeII}
        \[
          \scalebox{0.65}{\includegraphics{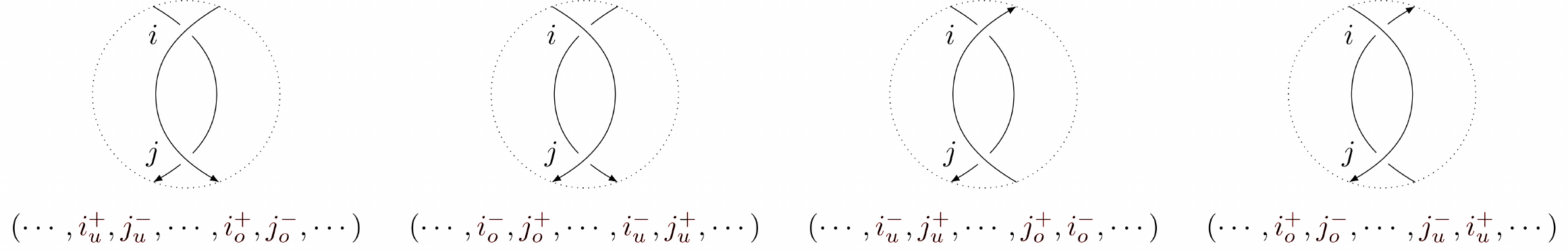}}
        \]
      \item Reidemeister III moves correspond to flipping three pairs. We have
        two cases, depending on the connectivity: \label{ax:reidGaussCodeIII}
        \[
          \scalebox{0.75}{\includegraphics{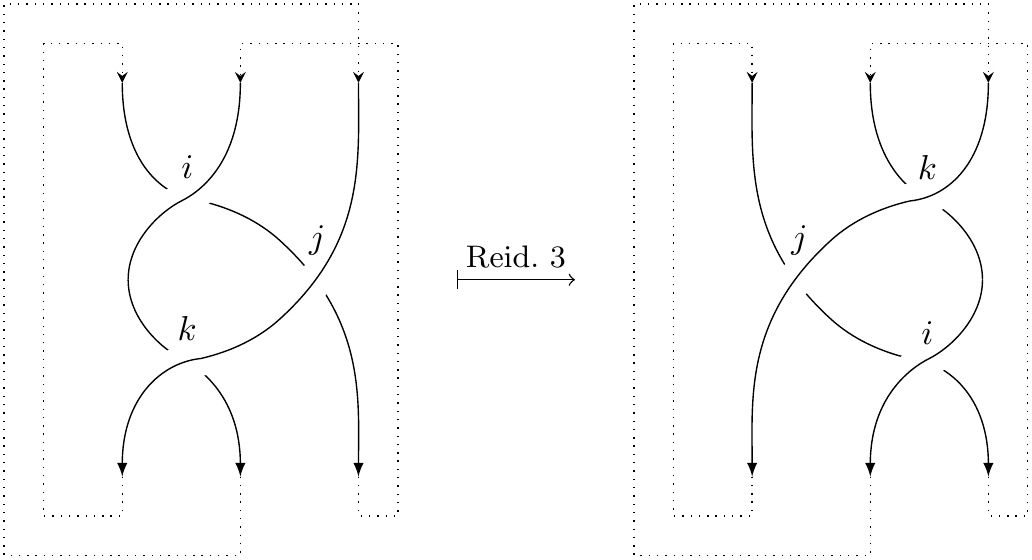}}\quad
          \scalebox{0.75}{\includegraphics{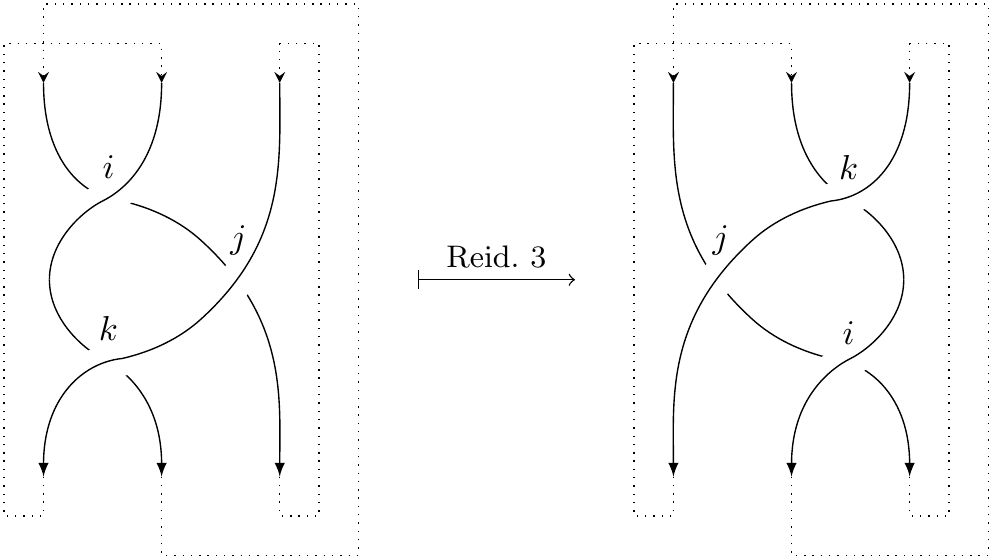}}
        \]
      \[
          \begin{tikzpicture}[scale=.8, every node/.style={scale=.8}]
            \def\myshift{5cm}
            \begin{scope}[xshift=-\myshift]
              \node (top) at (0,.5) {$(\cdots, i_u^+, j_u^+, \cdots, i_o^+,
                k_u^+, \cdots, j_o^+, k_o^+, \cdots)$};

              \node (bot) at (0,-.5) {$(\cdots, {  j_u^+},
                {  i_u^+}, \cdots, { k_u^+},
                { i_o^+}, \cdots, {  k_o^+}, {
                  j_o^+}, \cdots)$};

              \begin{scope}[xshift=-.3cm]
                \draw[-stealth] (-1.25, .25) -- (-1.75,-.25);
                \draw[white, line width=3pt] (-1.75, .25) -- (-1.25,-.25);
                \draw[-stealth] (-1.75, .25) -- (-1.25,-.25);

                \draw[-stealth] (.5, .25) -- (0,-.25);
                \draw[white, line width=3pt] (0, .25) -- (.5,-.25);
                \draw[-stealth] (0, .25) -- (.5,-.25);

                \draw[-stealth] (2.25, .25) -- (1.75,-.25);
                \draw[white, line width=3pt] (1.75, .25) -- (2.25,-.25);
                \draw[-stealth] (1.75, .25) -- (2.25,-.25);
              \end{scope}

            \end{scope}
            \begin{scope}[xshift=\myshift]
              \node (top) at (0,.5) {$(\cdots, i_u^+, j_u^+, \cdots, j_o^+,
                k_o^+, \cdots, i_o^+, k_u^+, \cdots)$};

              \node (bot) at (0,-.5) {$(\cdots, {  j_u^+},
                {  i_u^+}, \cdots, {  k_o^+}, {
                  j_o^+}, \cdots, {  k_u^+}, {  i_o^+},
                \cdots)$};

              \begin{scope}[xshift=-.3cm]
                \draw[-stealth] (-1.25, .25) -- (-1.75,-.25);
                \draw[white, line width=3pt] (-1.75, .25) -- (-1.25,-.25);
                \draw[-stealth] (-1.75, .25) -- (-1.25,-.25);

                \draw[-stealth] (.5, .25) -- (0,-.25);
                \draw[white, line width=3pt] (0, .25) -- (.5,-.25);
                \draw[-stealth] (0, .25) -- (.5,-.25);

                \draw[-stealth] (2.25, .25) -- (1.75,-.25);
                \draw[white, line width=3pt] (1.75, .25) -- (2.25,-.25);
                \draw[-stealth] (1.75, .25) -- (2.25,-.25);
              \end{scope}
            \end{scope}
          \end{tikzpicture}
        \]
    \end{itemize}

We note that Reidemeister moves do not change the parity of a crossing, so it
makes sense to regard parity as a fundamental property of a crossing.

\section{Biquandle Brackets}\label{BB}

Next, we recall a definition (see \cite{NOR,NO}).
  \begin{definition}
    Let $(X, \utr, \otr)$ be a biquandle and $R$ a commutative ring with
    identity. Let $w \in R^\times$, $\delta \in R$, and $A, B : X \times X \to
    R^\times$ such that the following hold (note, for notational
     compactness, we use $A_{x,y}$ as a shorthand for $A(x,y)$):
    \begin{itemize}
      \item For all $x \in X$, \label{ax:BQB1}
        \[
          -A_{x,x}^2 B^{-1}_{x,x} = w,
        \]
      \item For all $x,y \in X$, \label{ax:BQB2}
        \[
          -A_{x,y}^{-1} B_{x,y} - A_{x,y} B^{-1}_{x,y} = \delta,
        \]
      \item For all $x,y,z\in X$, we have \label{ax:BQB3}
        {
          \footnotesize
          \begin{align*}
            \btri{A}{A}{A}
            &= \btrI{A}{A}{A} \\[.5em]
            \btri{A}{B}{B}
            &= \btrI{B}{B}{A} \\[.5em]
            \btri{B}{B}{A}
            &= \btrI{A}{B}{B} \\[.5em]
            \btri{A}{B}{A}
            &= \phantom{+\delta} \btrI{A}{A}{B} + \btrI{B}{A}{A} \\
            &\pheq +\delta\btrI{B}{A}{B} + \btrI{B}{B}{B}
            \\[.5em]
            \btrI{A}{B}{A}
            &= \phantom{+\delta} \btri{A}{A}{B} + \btri{B}{A}{A}
            \\
            &\pheq + \delta\btri{B}{A}{B} + \btri{B}{B}{B}.
          \end{align*}
        }
    \end{itemize}
    Then we call $(A,B)$ an \emph{$X$-bracket over $R$}.
  \end{definition}

  The biquandle bracket axioms were chosen such that the \textit{state sum}
obtained by
  recursively simplifying a knot or link diagram $L$ according to the biquandle-colored skein relations
  \[
    \includegraphics{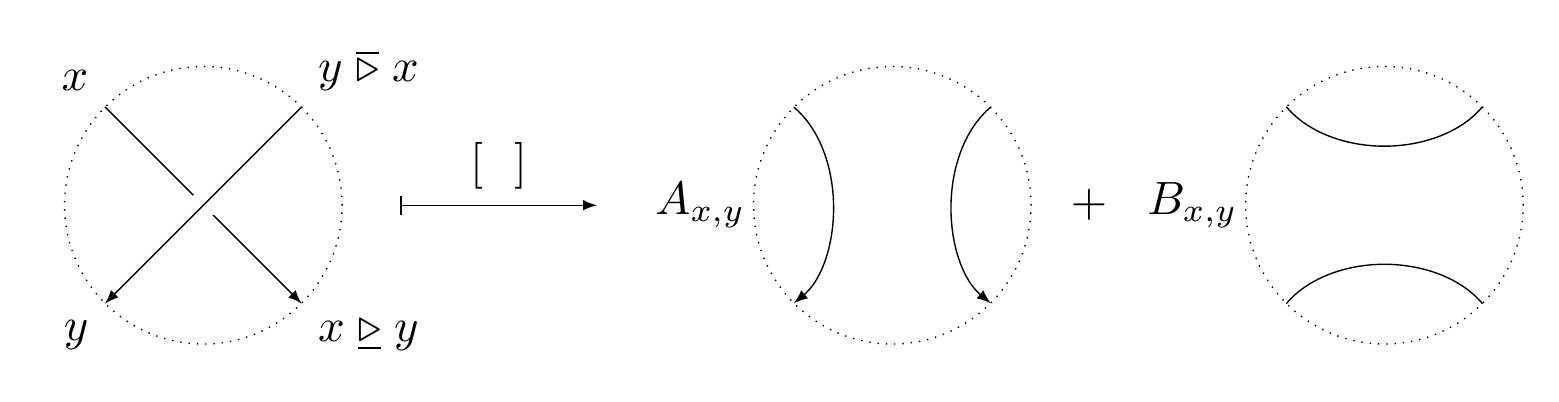}\]
\[
    \includegraphics{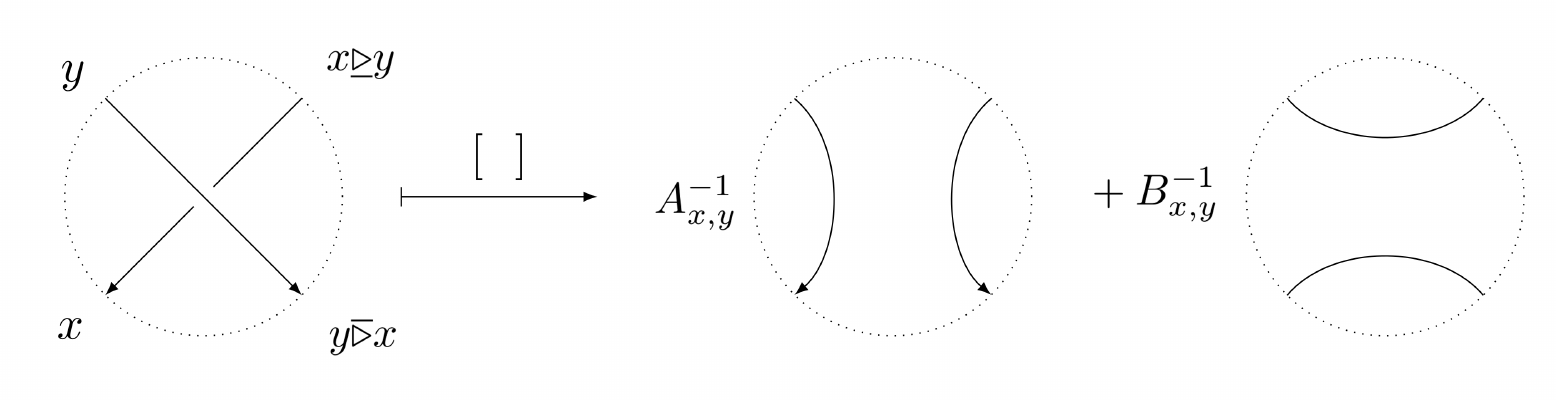}
  \]
then replacing each smoothed component with $\delta$ and multiplying
by the writhe correction factor $w^{n-p}$ (where $n,p$ are the numbers of
negative and positive crossings respectively) is invariant under
biquandle-colored Reidemeister moves. It follows that
the multiset of state-sum values over the set of biquandle colorings is an
invariant of oriented classical and virtual knots and links.

 Let us write $[D]$ for the state-sum value of a diagram $D$. Then it is straightforward to verify that we have
  \[
    \scalebox{0.75}{\includegraphics{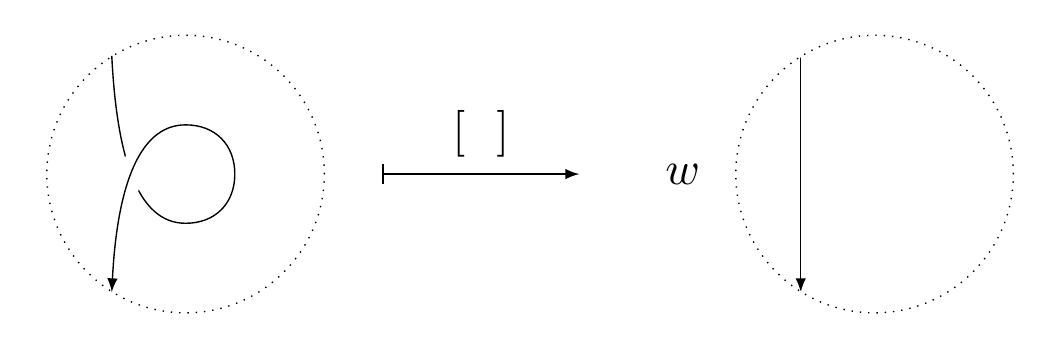}}
    \quad
    \scalebox{0.75}{\includegraphics{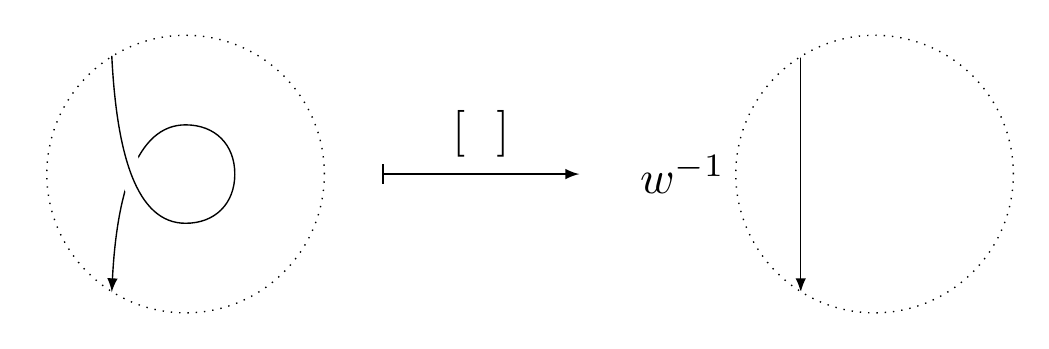}}
  \]
  as well as
  \[
    \includegraphics{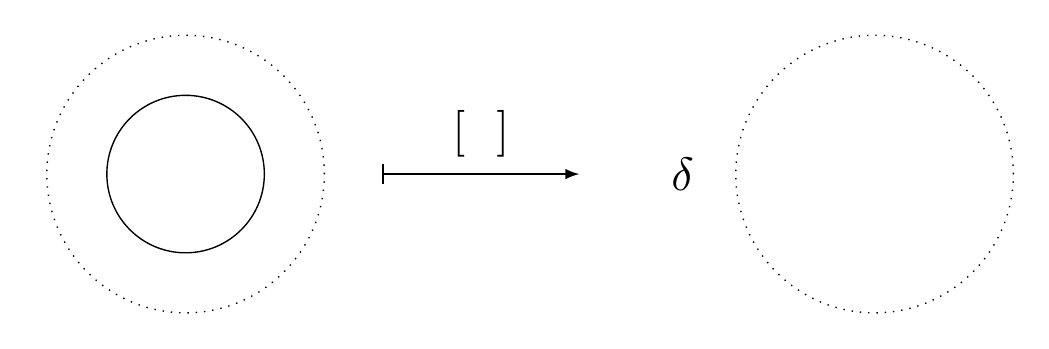}
  \]
  hence we can define the invariant as follows:
  \begin{definition}
 Let $L$ be an oriented knot or link diagram with a coloring by a finite biquandle $X$ and let $\beta$ be an $X$-bracket. Denote the number of positive crossings by $p$, the
        number of negative crossings by $n$, and the total number of crossings
        by $N$.
    \begin{itemize}
      \item A \emph{state} of $L$ is a choice of a type $A$ or type $B$
        smoothing at each crossing $j = 1, \ldots, N$ in $L$. We denote a state
        by $S \in \{A,B\}^{N}$, and denote the coefficient of the
        $j$\textsuperscript{th} smoothing in $S$ by $C^{S}_j$.
      \item For each state $S \in \{A,B\}^N$, let $m$ be the number of
        simple\footnote{Virtual crossings are permitted.} closed curves in $S$, and define
        \[
          C_{S} = \delta^m \prod_{j=1}^N C^{S}_j.
        \]
    \end{itemize}
    Then the \textit{state-sum} $\beta(L)$ given by
    \[
      \beta(L) = w^{n - p} \sum_{S \in \{A,B\}^N}  C_{S}
    \]
    is invariant under $X$-colored Reidemeister moves. See
\cite{NOR,NO} for more.
  \end{definition}

  Before moving on, it's worth noting that the smoothings defined above have two
  undesirable effects. First, both type $A$ and type $B$ smoothings can break
  colorings, as they concatenate incompatibly-colored semiarcs.
Second, type $B$ smoothings leave us with no
  well-defined orientation to our knots. This makes the idea of ``recursively''
  constructing the biquandle bracket invariant appear nonsensical.

  One possible solution is to simply ``perform all the smoothings
  simultaneously'' and obtain the state-sum form of the invariant. However, this
  loses one of the primary benefits of skein-relation-based invariants: namely,
  that partially-smoothed states can be simplified using the Reidemeister moves,
  thus reducing computations we need to perform.

  To address these problems, in \cite{NO} the second listed author and
collaborator reformulated biquandle brackets in terms of
oriented spatial graphs called \emph{trace diagrams}. The idea is to encode
  oriented links as graphs where the vertices correspond to crossings, and the
  directed edges correspond to oriented semiarcs. By extending biquandle
  colorings to these graphs, we obtain a straightforward way to formalize valid
  manipulations on partially-smoothed states.
  \begin{definition}
    A \emph{trace diagram} is an encoding of a partially smoothed
 oriented knot in a decorated directed graph as follows:
    \begin{itemize}
      \item Let $L$ be an oriented knot or link diagram with crossings
        $\{c_1, \ldots, c_n\}$. Since each crossing $c_i$ involves an
        overstrand and an understrand, we need to add \emph{two} corresponding
        vertices $u_i$ and $o_i$ to our digraph to determine which is
        which.\footnote{However in practice, we will simply draw $u_i$ and $o_i$
          on top of each other and use breaks in our edges to denote this
          instead.} Thus, let $V = \{u_1, o_1, \ldots, u_n,
          o_n\}$. \label{ax:TD1}
      \item For each oriented semiarc connecting crossings $(c_i, c_j)$ in $L$,
        add a directed edge from $u_i/o_i$ to $u_j/o_j$ as appropriate. If $L$
        is colored by some biquandle $X$, color all these edges
        accordingly. \label{ax:TD2}
      \item For each $c_i$, add a directed edge $(u_i, o_i)$ if $c_i$ is
        positive, and $(o_i, u_i)$ if $c_i$ is negative. We call these edges
        \emph{traces}, and will draw them in our diagrams with a dashed line.
        Traces corresponding to positive crossings are called \emph{positive}
        traces, and analogously in the negative case. \label{ax:TD3}
    \end{itemize}
     \end{definition}
  \begin{remark} In the case of classical knots and links, trace diagrams may be
visualized in $\mathbb{R}^3$ as follows:
    \[\scalebox{0.9}{
      \includegraphics{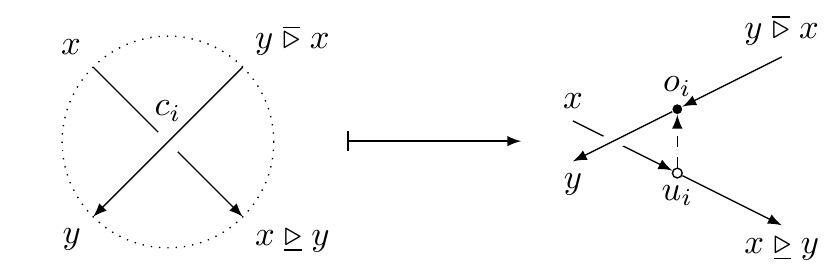}
      \hspace{1cm}
      \includegraphics{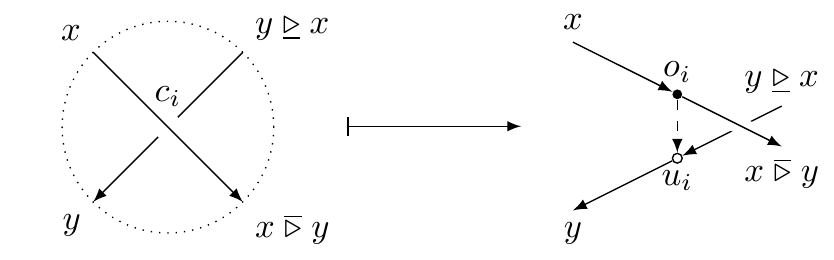}}
    \]
    The purpose of including the directed trace edges is to allow us to maintain biquandle colorings after we perform smoothings.
    Note that \ref{ax:TD3} has been chosen such that when drawing the $o_i$ on
    top of the $u_i$, positive traces point out of the page.
  \end{remark}
  \begin{example}
    The trefoil can be encoded in the following trace diagram:
    \[
      \includegraphics{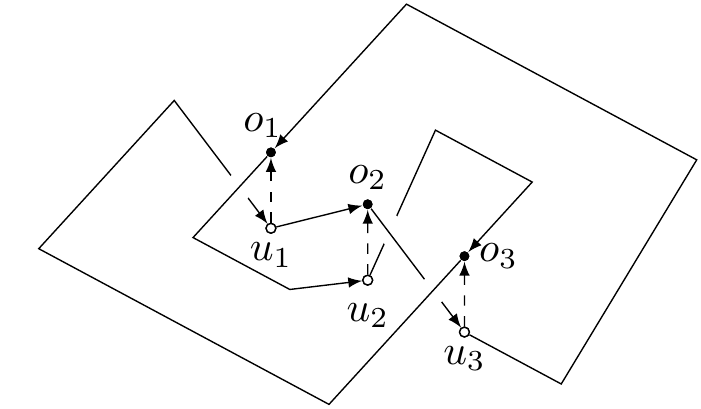}
    \]
    By drawing the $o_i, u_i$ on top of each other, the correspondence becomes
    easier to see:
    \[
      \includegraphics{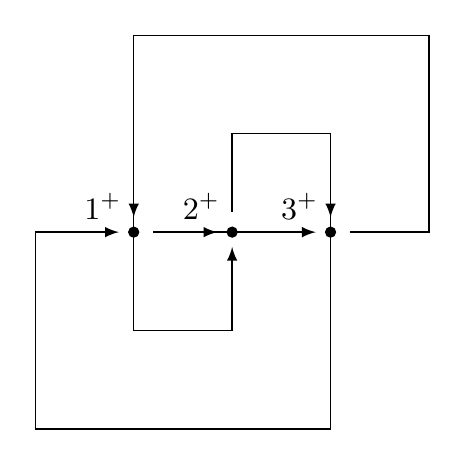}
    \]
  \end{example}
  Trace diagrams offer us a nice way to encode smoothings. We choose the
  following convention:
  \begin{definition}
    Let $X$ be a biquandle, and let $L$ be an $X$-colored knot/link with trace
    diagram $G = (V,A)$. Let $c_i$ be a crossing in $L$, and let $(\cdot, o_i),
    (o_i, \cdot), (\cdot, u_i), (u_i, \cdot)$, be the non-trace directed edges
    associated to $o_i$ and $u_i$. Then we define a \emph{smoothing} at $c_i$ by
    the rules

     \begin{center}
     \begin{tabular}{@{}ccc@{}}
       & \multicolumn{2}{c}{Smoothing Rule} \\ \hline 
       Sign of $c_i$ & Type $A$ & Type $B$ \\ \hline 
        $+$
        & \parbox{3cm}{\centering $(o_i, \cdot) \mapsto (u_i, \cdot)$ \\ $(u_i, \cdot) \mapsto (o_i, \cdot)$}
        & \parbox{3cm}{\centering $(\cdot, o_i) \mapsto (\cdot, u_i)$ \\ $(u_i, \cdot) \mapsto (o_i, \cdot)$}
        \\ \hline
        $-$
        & \parbox{3cm}{\centering $(\cdot, o_i) \mapsto (\cdot, u_i)$ \\ $(\cdot, u_i) \mapsto (\cdot, o_i)$}
        & \parbox{3cm}{\centering $(o_i, \cdot) \mapsto (u_i, \cdot)$ \\ $(\cdot, u_i) \mapsto (\cdot, o_i)$}
        \\\hline 
      \end{tabular}
    \end{center}
    which we can represent diagrammatically by
    \[
   \scalebox{0.9}{   \includegraphics{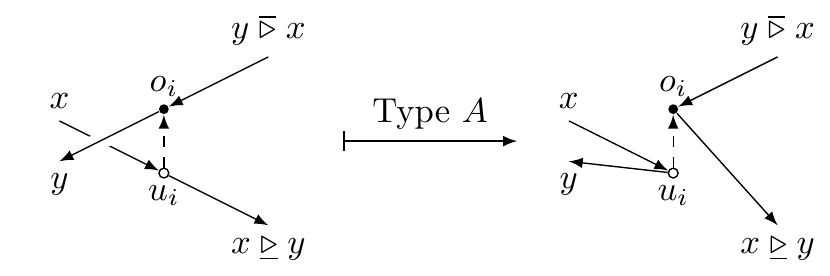}
      \hspace{1cm}
      \includegraphics{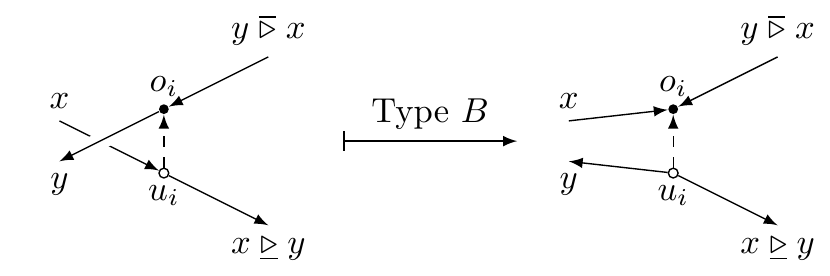} }
    \]
    for positive crossings, and
    \[
     \scalebox{0.9}{ \includegraphics{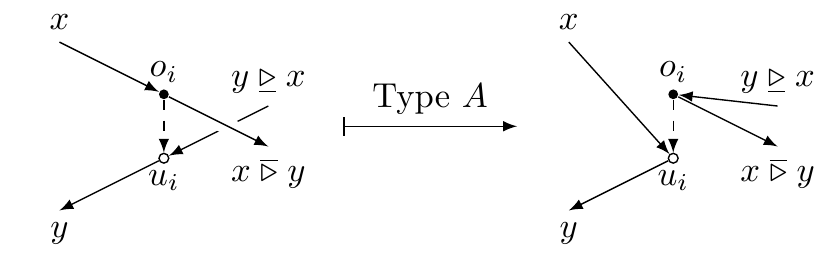}
      \hspace{1cm}
      \includegraphics{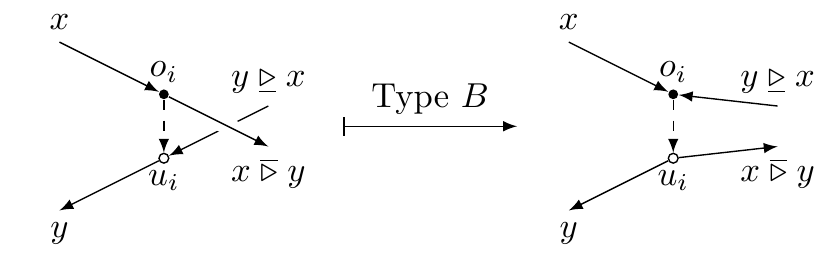}}
    \]
    for negative crossings. Flattening the $3$-D perspective, we see that these
    rules reflect the same smoothings employed in the biquandle bracket map:
    \[
    \scalebox{0.9}{  \includegraphics{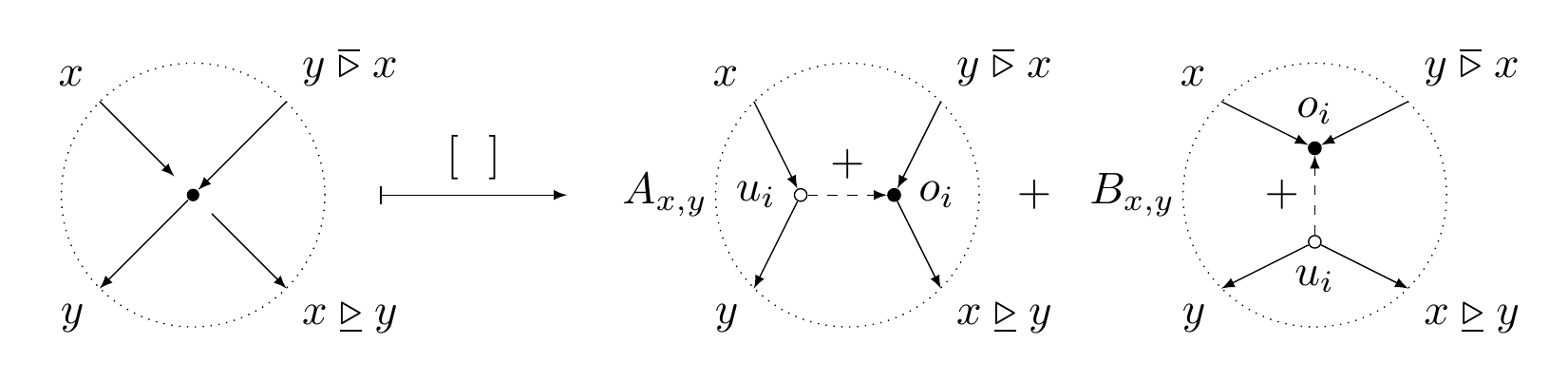}}\]
\[\scalebox{0.9}{
      \includegraphics{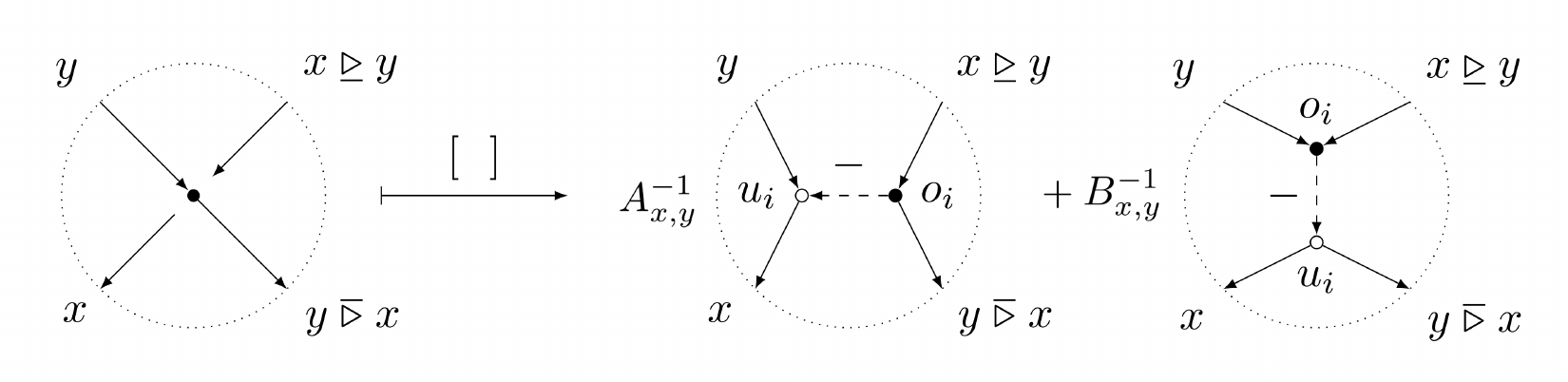}}
    \]
    Note that the rules we have chosen result in trace arrows pointing in
    positive $x,y$ for positive crossings, and negative $x,y$ for negative
    crossings.
  \end{definition}

  \begin{definition}
    Given a virtual knot or link $L$, we construct the trace diagram for $L$
    identically to the classical case, now decorating each vertex with parity
    information and denoting virtual crossings with circles.
  \end{definition}

\begin{remark}
In the case of virtual knots and links, the three-dimensional visualization
above should be understood locally within a neighborhood of each classical
crossing. Since the virtual Reidemeister moves do not change either the
biquandle colorings or connectivity in an oriented virtual link
diagram, it follows that starting with a biquandle-colored oriented virtual
link diagram we can apply the biquandle-bracket skein relations (either
with or without traces) to obtain a state-sum value $\beta$ which is invariant
under virtual Reidemeister moves.
\end{remark}

  \section{Kaestner Brackets}\label{KB}

We will now generalize biquandle brackets to the case of parity biquandles.
Note that the perhaps natural term for this construction, ``parity biquandle
bracket'', is in fact already in use in the literature with a different
meaning -- in \cite{IM}, the previously defined \textit{parity bracket}
construction from \cite{M}
was enhanced with biquandle colorings, whereas in our case we are enhancing
biquandle brackets by replacing
biquandles with parity biquandles. After briefly considering introducing
a distinction between ``(parity biquandle) brackets'' and
``parity (biquandle brackets)'', we opted instead to name the new structure
in honor of Aaron Kaestner, who introduced parity biquandles in \cite{KK}.

  \begin{definition}
    Let $X = (X, \utr^0, \otr^0, \utr^1, \otr^1)$ be a parity
    biquandle, and let $R$ be a commutative ring with identity. Let $\Ae, \Be,
    \Ao, \Bo : X \times X \to R^\times$. Then the collection $(X, \Ae, \Be, \Ao,
    \Bo)$ is a \emph{Kaestner bracket} if it satisfies the following
    conditions:
    \begin{itemize}
      \item $(X, \Ae, \Be)$ is a biquandle bracket,
      \item For all $x,y \in X$,
        \[
          \delta = -\pbA[1]{x,y} \cdot \pbB[1]{x,y}^{-1} - \pbA[1]{x,y}^{-1}
          \cdot \pbB[1]{x,y}
        \]
      \item For all $(\acol{a}, \bcol{b}, \ccol{c}) \in \{(1,1,0), (1,0,1),
          (0,1,1)\}$, we have the following (brace labels correspond to the trace
        diagrams in the following section): { \footnotesize
          \begin{align*}
            \overbrace{\tri{A}{A}{A}}^{\rm (i)}
            &= \overbrace{\trI{A}{A}{A}}^{\rm (I)} 
            \\
            \overbrace{\tri{A}{B}{B}}^{\rm (iv)}
            &= \overbrace{\trI{B}{B}{A}}^{\rm (VII)} 
            \\
            \overbrace{\tri{B}{B}{A}}^{\rm (vii)}
            &= \overbrace{\trI{A}{B}{B}}^{\rm (IV)} 
            \\
            \overbrace{\tri{A}{B}{A}}^{\rm (iii)}
            &=
              \overbrace{\trI{A}{A}{B} \phantom{\delta} }^{\rm (II)} +
              \overbrace{\trI{B}{A}{A}}^{\rm (V)} 
            \\
            & +
              \underbrace{\delta\trI{B}{A}{B}}_{\rm (VI)} +
              \underbrace{\trI{B}{B}{B}}_{\rm (VIII)}
            \\
            \overbrace{\trI{A}{B}{A}}^{\rm (III)}
            &=
              \overbrace{\tri{A}{A}{B} \phantom{\delta} }^{\rm (ii)} +
              \overbrace{\tri{B}{A}{A}}^{\rm (v)} 
            \\
            & +
              \underbrace{\delta\tri{B}{A}{B}}_{\rm (vi)} +
              \underbrace{\tri{B}{B}{B}}_{\rm (vii)}\nonumber
            \\
          \end{align*}
        }
      \end{itemize}
      \vspace{-1cm}
  \end{definition}
  This definition is motivated by the following skein relations
\[\includegraphics{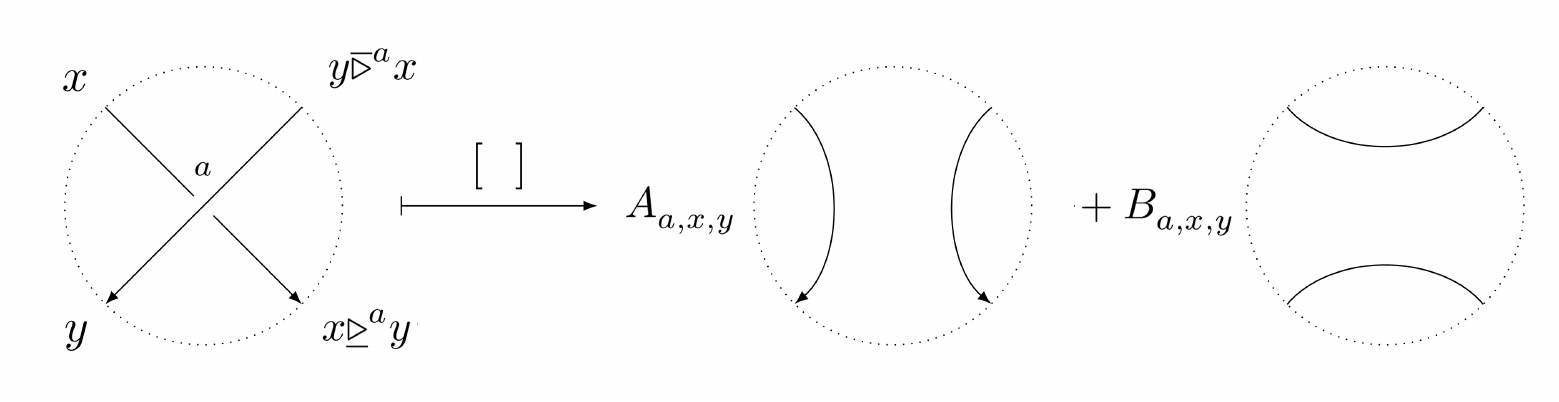}\]
\[\includegraphics{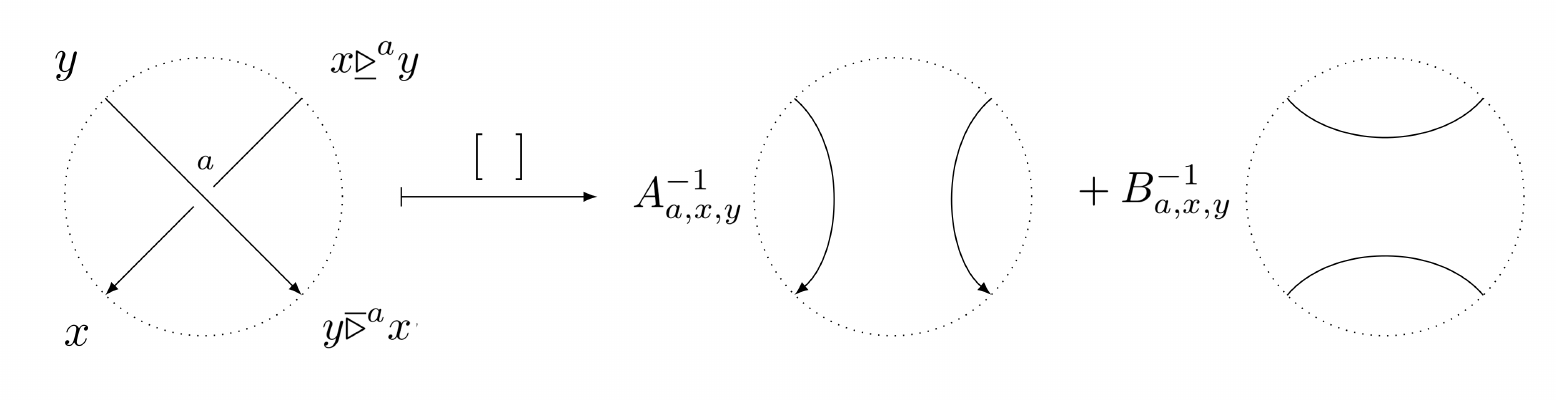}\]
with the same values of $\delta$ and $w$ as before. Then we obtain the
conditions above by considering
  trace diagram reductions of both sides of the positive Reidemeister III move.
  Applying the smoothing rules, the LHS gives
\[ \begin{array}{cccccccc}
   \scalebox{1.35}{\includegraphics{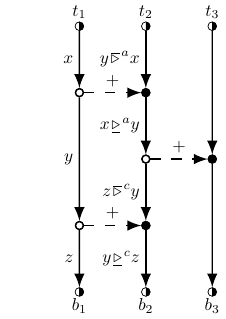}} &
   \scalebox{1.35}{\includegraphics{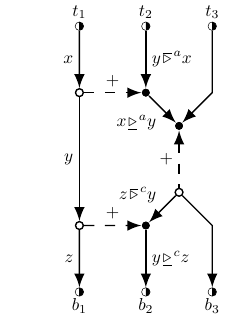}} &
   \scalebox{1.35}{\includegraphics{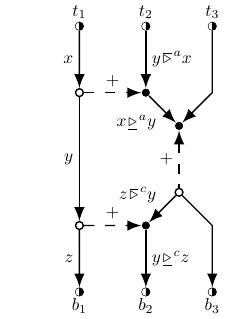}} &
   \scalebox{1.35}{\includegraphics{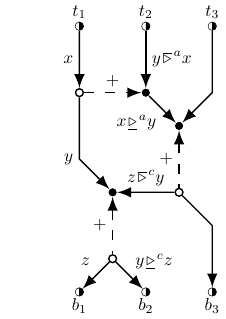}} \\
\mathrm{(i)} & \mathrm{(ii)} & \mathrm{(iii)} &\mathrm{(iv)}\\
   \scalebox{1.35}{\includegraphics{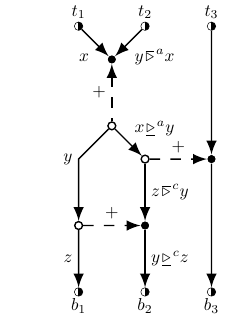}} &
   \scalebox{1.35}{\includegraphics{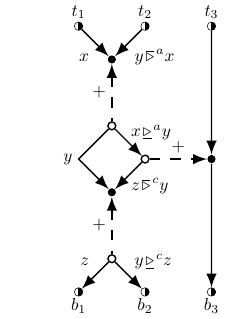}} &
   \scalebox{1.35}{\includegraphics{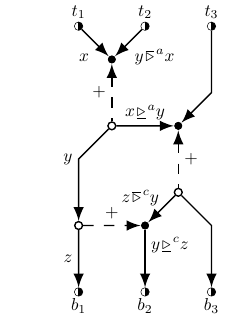}} &
   \scalebox{1.35}{\includegraphics{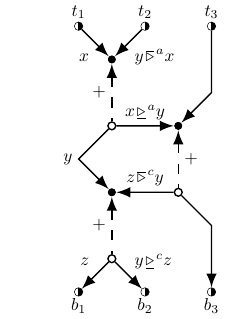}}\\
\mathrm{(v)} &\mathrm{(vi)} &\mathrm{(vii)} &\mathrm{(viii)}
\end{array}
\]
    while the RHS gives
\[\begin{array}{cccccccc}
   \scalebox{1.2}{\includegraphics{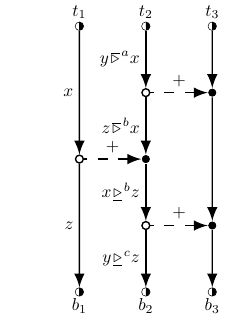}} &
   \scalebox{1.2}{\includegraphics{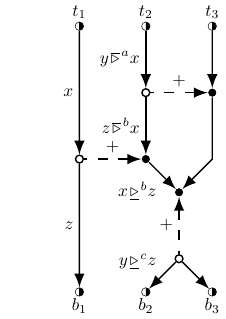}} &
   \scalebox{1.2}{\includegraphics{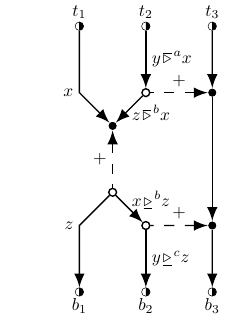}} &
   \scalebox{1.2}{\includegraphics{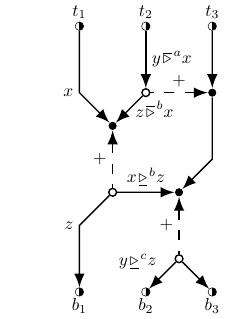}} \\
\mathrm{(I)} & \mathrm{(II)} & \mathrm{(III)} &\mathrm{(IV)}\\
   \scalebox{1.2}{\includegraphics{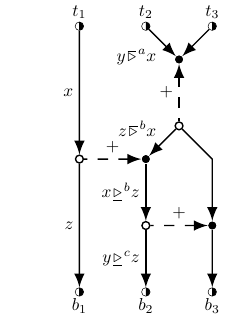}} &
   \scalebox{1.2}{\includegraphics{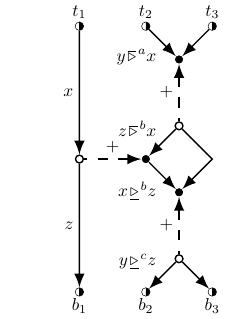}} &
   \scalebox{1.2}{\includegraphics{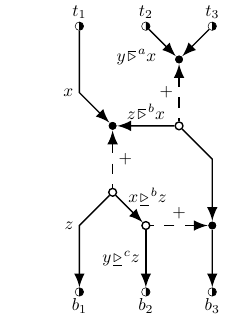}} &
   \scalebox{1.2}{\includegraphics{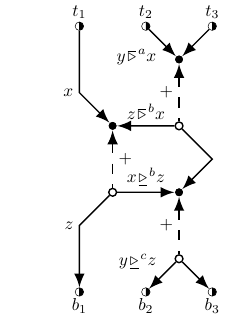}} \\
\mathrm{(V)} &\mathrm{(VI)} &\mathrm{(VII)} &\mathrm{(VIII)} \end{array}
\]

\begin{example}\label{ex:kb1}
Let $X$ be the parity biquandle in Example \ref{ex:pb1} and let
$R=\mathbb{Z}_5$. The reader can verify that the the coefficients in the tables
\[\begin{array}{r|rrr}
A_{0} & 1 & 2 & 3 \\ \hline 
   1 & 1 & 4 & 4\\ 
   2 & 4 & 1 & 4\\ 
   3 & 4 & 4 & 1 \end{array}\quad
\begin{array}{r|rrr}
B_{0} & 1 & 2 & 3 \\ \hline
   1 & 4 & 1 & 1\\ 
   2 & 1 & 4 & 1\\ 
   3 & 1 & 1 & 4 \end{array}\quad
\begin{array}{r|rrr}
A_{1} & 1 & 2 & 3 \\ \hline
   1 & 1 & 3 & 1\\ 
   2 & 4 & 2 & 1\\ 
   3 & 4 & 3 & 4 \end{array}\quad
\begin{array}{r|rrr}
B_{1} & 1 & 2 & 3 \\ \hline
   1 & 4 & 2 & 4\\ 
   2 & 1 & 3 & 4\\ 
   3 & 1 & 2 & 1 \end{array}\]
define a Kaestner bracket over $R$ with $w=-(1)^2(4^{-1})=1$ and
$\delta=-(1^{-1})(4)-4^{-1}(1)=1+1=2$.
This data encodes 54 skein relations
depending on the parity, crossing sign and colors at a crossing, including
for example
\[\includegraphics{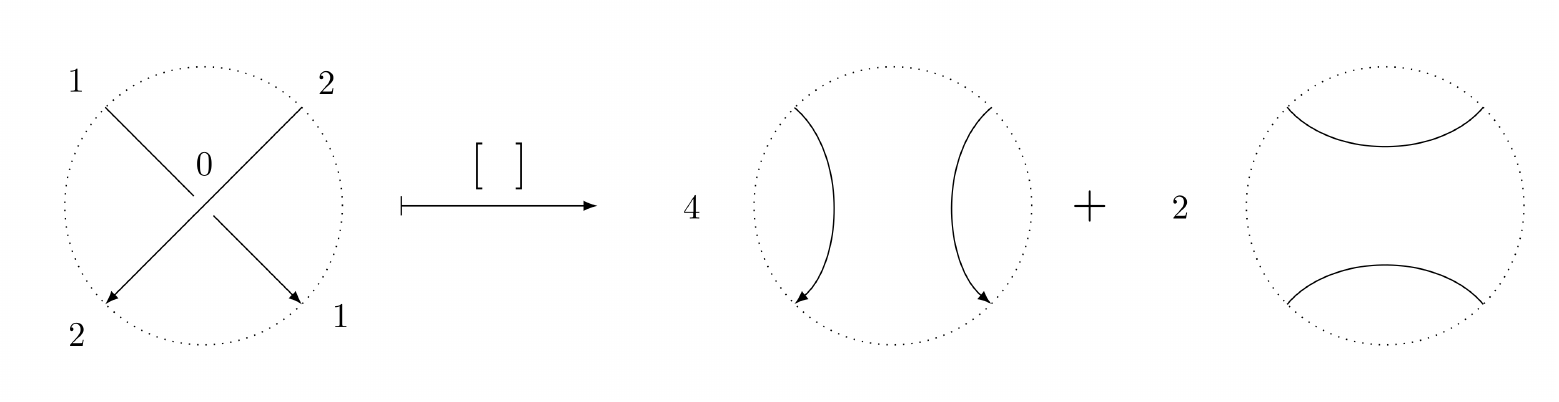}\]
and
\[\includegraphics{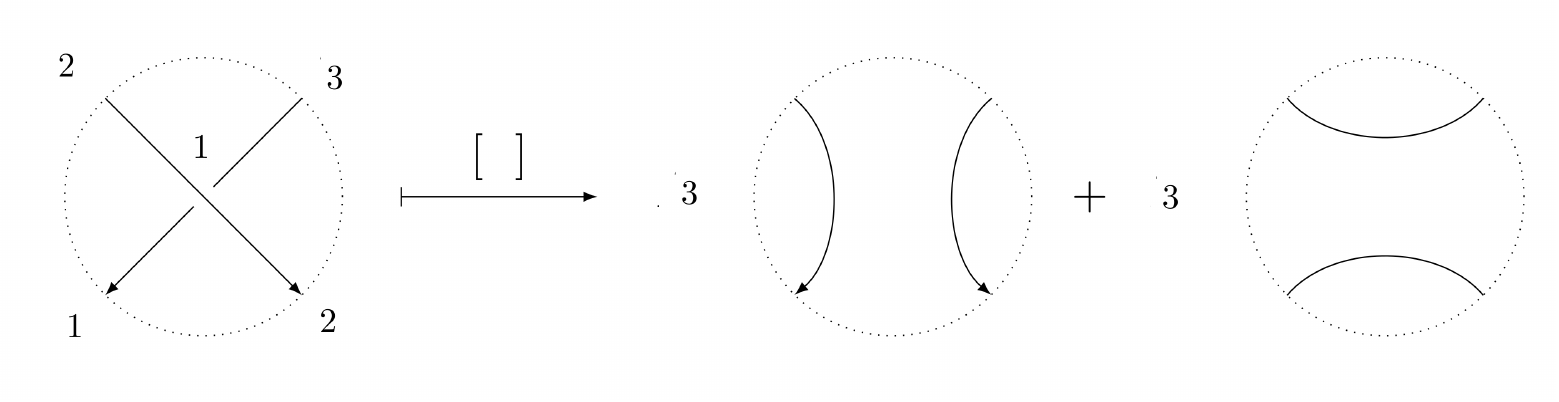}\]
\end{example}

We can now state our main definition.

\begin{definition}
Let $X$ be a parity biquandle, $R$ a commutative ring with identity, and
$\beta$ a Kaestner bracket on $X$ with coefficients in $R$. Then for any
oriented virtual link $L$, we define the \textit{Kaestner bracket polynomial}
of $L$ to be the sum of contributions of a formal variable $u$ to the power
of the state sum $\beta(L_f)$ over the set of $X$-colorings $L_f$ of $L$,
i.e.,
\[\Phi_X^{\beta}(L)=\sum_{L_f\in\mathcal{C}(L,X)} u^{\beta(L_f)}.\]
\end{definition}

By construction, we have the following:
\begin{proposition}
For any finite parity biquandle $X$ and Kaestner bracket $\beta$,
$\Phi_X^{\beta}(L)$ is a invariant of oriented virtual links.
\end{proposition}

\begin{example}\label{ex:kb2}
The virtual trefoil knot $2.1$ has three colorings by the parity biquandle
in Example \ref{ex:pb1} as shown.
\[\includegraphics{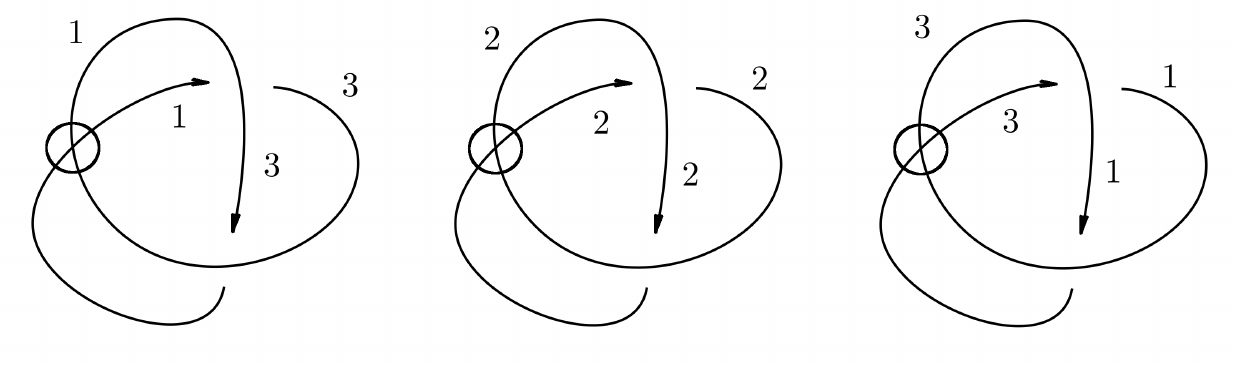}\]
Let us compute the state-sum value $\beta$ for the first one.
\[\includegraphics{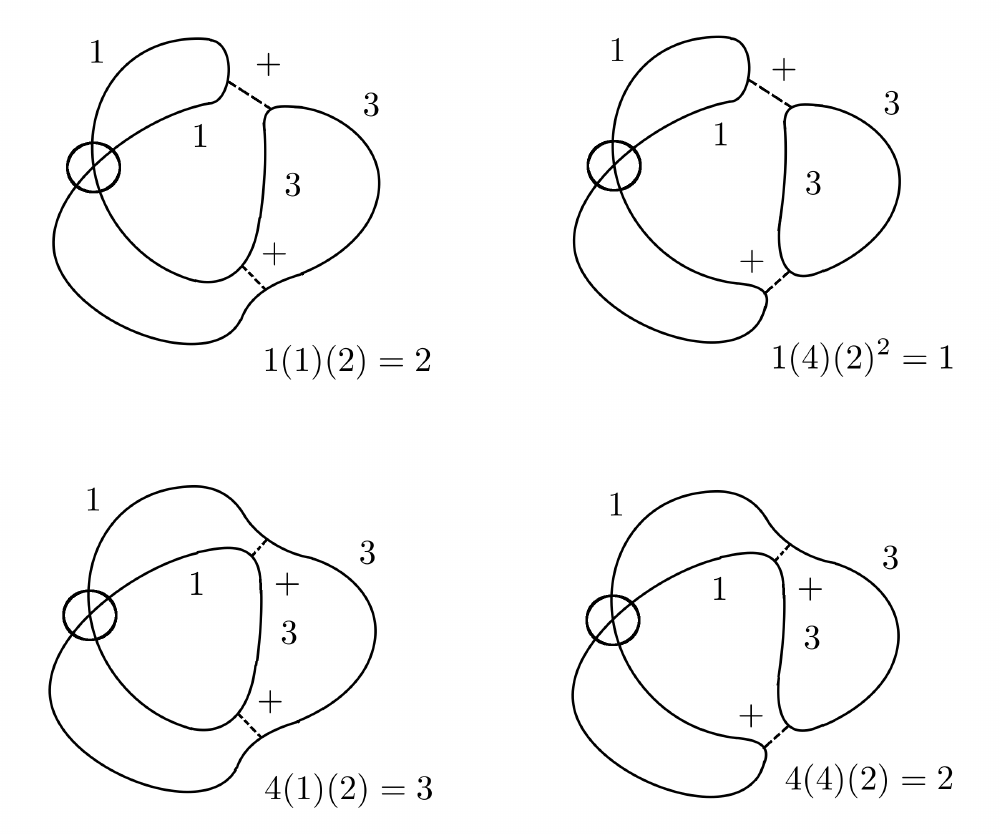}\]
The four states contribute $2+1+3+2=3$ for a state-sum contribution of
$u^3$ to the invariant value.
Repeating with other two colorings, we obtain $\Phi_X^{\beta}(2.1)=3u^3$. This
result distinguishes $2.1$ from the unknot, which has invariant
value $\Phi_X^{\beta}(\mathrm{unknot})=3u^2$.
\end{example}

\begin{example}\label{ex:kb3}
We computed $\Phi_X^{\beta}$ using the parity biquandle and Kaestner bracket
from Example \ref{ex:kb1} for the virtual knots with up to four classical 
crossings in the knot atlas \cite{KA}.
The results are in the table. Note that the counting invariant alone does not
distinguish any of the virtual knots on the table, while this particular
Kaestner bracket sorts them into two classes.
\[\begin{array}{r|l}
\Phi_X^{\beta}(K) & K \\\hline
3u^2 & 3.1, 3.5, 3.6, 3.7, 4.1, 4.2, 4.3, 4.6, 4.7, 4.8, 4.10, 4.12, 4.13, 4.16, 4.17, 4.19, 4.21, 4.23, 4.24, 4.25, 4.26, \\ 
& 4.31, 4.32, 4.35, 4.36, 4.41, 4.42, 4.43.4.46, 4.47, 4.50, 4.51, 4.53, 4.55, 4.56, 4.57, 4.58, 4.59, 4.65, 4.66, \\ &
4.67, 4.68, 4.70, 4.71, 4.72, 4.73, 4.75, 4.76, 4.77, 4.79, 4.80, 4.85, 4.86, 4.89, 4.90, 4.91, 4.93, 4.96, 4.97, \\ &
4.98, 4.99, 4.100, 4.102, 4.103, 4.105, 4.106, 4.107, 4.108 \\ \hline
3u^3 & 2.1, 3.2, 3.3, 3.4, 4.4, 4.5, 4.9, 4.11, 4.14, 4.15, 4.18, 4.20, 4.22, 4.27, 4.28, 4.29, 4.30, 4.33, 4.34, \\ &
4.37, 4.38, 4.39, 4.40, 4.44, 4.45, 4.48, 4.49, 4.52, 4.54, 4.60, 4.61, 4.62, 4.63, 4.64, 4.69, 4.74, 4.78, 4.81, \\ &
4.82, 4.83, 4.84, 4.87, 4.88, 4.92, 4.94, 4.95, 4.101, 4.104\\
\end{array}\]
In particular, treating all crossings as even and evaluating the resulting
classical biquandle bracket in this example yields the value $3u^2$ for all of
the virtual knots on the table; hence this example shows that Kaestner brackets
in general define stronger invariants than their classical biquandle 
bracket counterparts. Moreover, since the parity biquandle counting invariant
value is $3$ for all virtual knots on the table, this example shows that
Kaestner brackets are stronger invariants than the unenhanced parity 
biquandle counting invariants. 
\end{example}

  \section{Questions}\label{Q}

  We conclude in this section with some interesting questions and possible 
directions for future research.
    \begin{itemize}
    \item Is the current list of axioms minimal? It seems likely that the
      answer is ``no.'' Finding a simpler set of axioms could improve the
      performance of computer search programs by allowing fewer condition
      checks.
    \item Limiting our search to brackets with coefficients in finite rings
allows for computer searches; we are very interested in finding brackets
with coefficients in infinite rings ($\mathbb{Z}$, polynomials over infinite
fields, etc.)
    \item What ways are there to modify the value $\delta$ of a state component
      to depend on the biquandle coloring as well? 
    \item What properties of Kaestner brackets can be used to detect non-classicality or almost-classicality? 
    \end{itemize}

  \bibliography{fk-sn2-rev1}{}
  \bibliographystyle{abbrv}

  \noindent
  \textsc{Department of Mathematics \\
  Harvey Mudd College\\
  301 Platt Boulevard \\
  Claremont, CA 91711
  }

  \bigskip

  \noindent
  \textsc{Department of Mathematical Sciences \\
  Claremont McKenna College \\
  850 Columbia Ave. \\
  Claremont, CA 91711}

\end{document}